\documentclass{article}

\newif\ifpdf
\ifx\pdfoutput\undefined \pdffalse \else \pdfoutput=1 \pdftrue \fi

\usepackage{amssymb}
\usepackage{amsmath}                    
\usepackage{amsthm}

\ifpdf
\usepackage[pdftex]{graphicx}
\usepackage[pdftex]{hyperref}
\else
\usepackage{graphicx}
\fi

\theoremstyle{plain}
\newtheorem*{lebrun}{Theorem (C. LeBrun)}
\newtheorem*{chen}{Theorem (Z. Chen)}
\newtheorem*{hitchin-thorpe}{Theorem (N. Hitchin, J. Thorpe)}
\newtheorem*{gromov}{Theorem (M. Gromov)}
\newtheorem*{theoremA}{Theorem A}
\newtheorem*{theoremB}{Theorem B}
\newtheorem{theorem}{Theorem}
\newtheorem{lemma}[theorem]{Lemma}
\newtheorem{proposition}[theorem]{Proposition}
\newtheorem{corollary}[theorem]{Corollary}

\theoremstyle{definition}
\newtheorem{definition}{Definition}

\theoremstyle{remark}
\newtheorem*{remark}{Remark}

\newcommand{\mfld}[1]{${#1}$-ma\-ni\-fold}
\newcommand{\CP}[1]{\mathbb{CP}^{#1}}

\newcommand{\F}{\mathbb{F}}
\newcommand{\N}{\mathbb{N}}
\newcommand{\R}{\mathbb{R}}
\newcommand{\Z}{\mathbb{Z}}

\DeclareMathOperator{\spin}{\ensuremath{Spin}}

\DeclareMathOperator{\supp}{\ensuremath{supp}}
\DeclareMathOperator{\hol}{\ensuremath{hol}}
\DeclareMathOperator{\inj}{\ensuremath{inj}}

\allowdisplaybreaks[1]

\begin{document}

\title{Seiberg-Witten invariants of non-simple type and Einstein metrics}
\author{Heberto del Rio Guerra\thanks{The author acknowledges support from grant 28491-E of CONACyT}\\Centro de Investigaci\'on en Matem\'aticas, A.C. (CIMAT)\\Guanajuato, M\'exico}
\date{\today}

\maketitle

\begin{abstract} 
We construct examples of four dimensional manifolds with 
$\spin^c$-struc\-tures, whose moduli spaces of solutions to the Seiberg-Witten equations, represent a non-trivial bordism class of 
positive dimension, {\it i.e.} the $\spin^c$-structures are not induced by almost complex structures. 
As an application, we show the existence of infinitely many non-homeomorphic compact oriented 
$4$-manifolds with free fundamental group and predetermined Euler characteristic and signature that do not carry Einstein metrics (see \cite{sambusetti}).
\end{abstract}

\section{Introduction}

A smooth Riemannian manifold $(M,g)$ is said to be \emph{Einstein} if its Ricci curvature tensor $r$ is a multiple of the metric \textit{i.e.}
\begin{equation*}
r=\lambda g.
\end{equation*}
Not every smooth compact oriented \mfld{4} admits such a metric. A well known obstruction is given by the following result due to N. Hitchin and J. Thorpe (see \cite{besse2}). \emph{ If $M$ is a compact oriented \mfld{4} and $e(M)<\frac{3}{2}|\sigma(M)|$ then $M$ does not admit an Einstein metric, where $e$ and $\sigma$ respectively denote the Euler characteristic and the signature.} The Gauss-Bonnet-like formula
\begin{equation*}
2e(M)\pm3\sigma(M)=\frac{1}{4\pi^2}\int_M\left(\frac{s^2}{24}-\frac{|r_0|^2}{2}+2|W_\pm|^2\right)d\mu,
\end{equation*}
implies Hitchin-Thorpe's inequality because Einstein metrics are characterized by the vanishing of $r_0$, and this is the only negative term in the above integrand. Here $s,\ r_0,\ W_+,\ W_-$ respectively denote the scalar, trace-free Ricci, self-dual Weyl, and anti-self-dual Weyl curvature tensors of a Riemannian metric. 

As C. LeBrun showed in \cite{lebrun} this result can be improved using careful estimates on the $L^2$-norm of the scalar curvature tensor $s$ and the $L^2$-norm of the self-dual part of the Weyl tensor $W_+$ arising from the Seiberg-Witten equations if, for example, the smooth \mfld{4} $M$, admits a symplectic form. To obtain these estimates C. LeBrun used that such an $M$ admits irreducible solutions to the Seiberg-Witten equations for every metric $g$ rather than actually using the fact that $M$  has non-trivial Seiberg-Witten invariant. Our main result is
\begin{theoremA}
Let $(M,\mathfrak{c})$ be a smooth compact K\"ahler surface with a $\spin^c$-structure
$\mathfrak{c}$. There is a canonical $\spin^c$ structure in the connected sum manifold
$M\#(S^1\times S^3)$ which we will denote by $\mathfrak{c}_{0,1}$. Moreover
$d(\mathfrak{c}_{0,1})=d(\mathfrak{c})+1$. If $\mathfrak{c}$ is a non-trivial SW-class
for $M$ then $\mathfrak{c}_{0,1}$ is a $B$-class for the connected sum $M\#(S^1\times
S^3)$.
\end{theoremA}

The equality $d(\mathfrak{c}_{0,1})=d(\mathfrak{c})+1$ implies that $\mathfrak{c}_{0,1}$ is not induced by an almost complex structure, and the statement \emph{$\mathfrak{c}_{0,1}$ is a $B$-class} implies that there exist irreducible solutions to the Seiberg-Witten equations for every Riemannian metric. The technique that we have used to produce these $\spin^c$-structures does not rely on the well-known gluing-argument (compare with \cite{ozsvath-szabo}). 

The main application of our result is
\begin{theoremB}
For each admissible pair $(m,n)$ there exist an infinite number of non-homeomorphic compact oriented $4$-manifolds which have Euler characteristic $m$, signature $n$, with free fundamental group and which do not admit an Einstein metric.
\end{theoremB}
Similar examples but with very complicated fundamental group have been obtained by A. Sambusetti
\cite{sambusetti} using connected sums with real or complex hyperbolic $4$-manifolds.

I would like to thank Prof. C. LeBrun for all the useful comments, and the time he spent reading previous versions of this manuscript.

\section{SW-Moduli Space}

\begin{definition}\label{defsw}
Let $(M,\mathfrak{c})$ be a smooth compact oriented \mfld{4} with
a $\spin^c$-structure $\mathfrak{c}$. Let
$L_\mathfrak{c}=\det(\mathfrak{c})$ be the determinant line bundle
associated to $\mathfrak{c}$. Fix a Riemannian metric $g$ on $M$.
The configuration space $\mathcal{C}(\mathfrak{c})$ consist of
pairs $(A,\phi)$ , where $A$ is an $U(1)$-connection on
$L_\mathfrak{c}$ and
$\phi\in\mathcal{C}^\infty(S^+(\mathfrak{c}))$ is a self-dual
spinor. We say that $(A,\phi)$ satisfy the \mbox{Seiberg-Witten}
equations (SW-equations) if and only if
\begin{align*}
D_A\phi&=0\\ F_A^+&=q(\phi),
\end{align*}
where $q(\phi)=\phi\otimes\phi^*-\frac{|\phi|^2}{2}\text{Id}$.
\end{definition}
\begin{remark}
$D_A$ is the associated Dirac operator of the $\spin^c$-bundle,
and $F_A^+$ is the self-dual part of the curvature associated to
the connection $A$, thought of as an endomorphism of the self-dual spinors.
\end{remark}

\begin{definition}
We say that an element $(A,\phi)$ is irreducible if
$\phi\not\equiv 0$, otherwise it is reducible. We denote by
$\mathcal{C}^*(\mathfrak{c})$ the open subset of irreducible
configurations, by $\mathcal{G}(\mathfrak{c})=\{\sigma:M\to S^1\}$
the gauge group, and by
$\mathcal{B}^*(\mathfrak{c})=\mathcal{C}^*(\mathfrak{c})/\mathcal{G}(\mathfrak{c})$
the open subset of irreducible equivalence classes.
\end{definition}

The naive definition of the \mbox{Seiberg-Witten} moduli space would be:
\begin{equation*}
\mathcal{M}_g(\mathfrak{c})=\{(A,\phi)\in\mathcal{C}(\mathfrak{c})|\ D_A\phi=0,\
F_A^+=q(\phi)\}/\mathcal{G}(\mathfrak{c}),
\end{equation*}
but in order to use the usual analytical tools, one has to extend
the $\mathcal{C}^\infty$ objects to appropriate Sobolev spaces.
From now on we extend the configuration space
$\mathcal{A}(\mathfrak{c})$ and the gauge group
$\mathcal{G}(\mathfrak{c})$ by requiring $A$ and $\phi$ to be in
$L^2_2$ and $\sigma$ to be in $L^2_3$. The SW-equations and the
gauge actions make sense in this context also and we define:

\begin{definition} The \mbox{Seiberg-Witten} moduli space is:
\begin{equation*}
\mathcal{M}_g(\mathfrak{c})=\{(A,\phi)\in\mathcal{C}(\mathfrak{c})|\ D_A\phi=0,\
F_A^+=q(\phi)\}/\mathcal{G}(\mathfrak{c}),
\end{equation*}
where $\mathcal{A}(\mathfrak{c})$ and $\mathcal{G}(\mathfrak{c})$ are the extended
configuration space and gauge group. The formal dimension (computed using the Atiyah-Singer index theorem) of this moduli space is
\begin{equation*}
d(\mathfrak{c})=\frac{c_1^2(\mathfrak{c})-(2e(M)+3\sigma(M))}{4}.
\end{equation*}
\end{definition}
In general there is no reason to expect that the moduli space form a smooth manifold. The
best we can hope for is that \emph{generically} it does. The next Theorem guarantees that
this is the case. For the proof see \cite{morgan:sw}.

\begin{theorem}
Suppose that $b_2^+>0$. Fix a metric $g$ on $M$. Then for a generic $\mathcal{C}^\infty$
self-dual $2$-form $h$ on $M$ the following holds. For any $\spin^c$-structure
$\mathfrak{c}$ on $M$ the moduli space
$\mathcal{M}_g(\mathfrak{c},h)\subset\mathcal{B}(\mathfrak{c})$ of gauge equivalence
classes of pairs $[A,\phi]$ which are solutions to the perturbed SW-equations
\begin{align*}
D_A\phi&=0\\ F_A^+-q(\phi)&=ih
\end{align*}
form a smooth compact submanifold of $\mathcal{B}^*(\mathfrak{c})$ of dimension
$d(\mathfrak{c})$.
\end{theorem}
Also in \cite{morgan:sw} it is shown that if $b_2^+>1$ then the bordism class of $\mathcal{M}_g(\mathfrak{c},h)$ is an invariant of the smooth structure of $M$ and the $\spin^c$-structure $\mathfrak{c}$ on $M$. We will denote by  $\mathcal{M}(\mathfrak{c})$ this bordism class.

\begin{proposition}\label{cont}
Consider a fixed $U(1)$-connection $A$ on $L_\mathfrak{c}$. Let $[A_i,\phi_i]$ be
solutions to the SW-equations, and let $(A_i,\phi_i)$ be the unique representatives such
that $A_i-A$ is co-closed (gauge fixing condition, see \cite{morgan:sw}), for $i=1,2$. If
$\phi_1=\phi_2$ then $A_1=A_2$.
\end{proposition}
\begin{proof}
The first thing to notice is that $A_2=A_1+\theta$, where $\theta$ is a co-closed
$1$-form. Since $(A_1,\phi_1)$ and $(A_2,\phi_2)$ are solutions to the SW-equations we
have
\begin{align*}
F_{A_1}^+&=q(\phi_1)\\ &=q(\phi_2)\\ &=F_{A_2}^+.
\end{align*}
Therefore
\begin{align*}
F_{A_2}^+-F_{A_1}^+=0&\Leftrightarrow (d\theta)^+=0\\ &\Leftrightarrow
*d\theta=-d\theta\\ &\Rightarrow d*d\theta=-dd\theta=0\\ &\Leftrightarrow *d*d\theta=0\\
&\Leftrightarrow \delta d\theta=0.
\end{align*}
This last statement and the fact $\delta\theta=0$ implies that
\begin{align*}
\Delta\theta&=d\delta\theta+\delta d\theta\\ &=\delta d\theta\\ &=0.
\end{align*}
Since $(A_i,\phi_i)$ $i=1,2$ are solutions to the \mbox{Seiberg-Witten} equations we have
\begin{align*}
0&=D_{A_2}\phi_2\\ &=D_{A_1+\theta}\phi_1\\ &=D_{A_1}\phi_1+\theta\cdot\phi_1\\
&=\theta\cdot\phi_1,
\end{align*}
\emph{multiplying} by $\theta$ both sides of the equality we get that
$|\theta|^2\phi_1=0$. Taking the point-wise norm we will have $|\theta|^2|\phi_1|=0$. If
we denote by $Z_{|\theta|^2}$ and $Z_{|\phi_1|}$ the set of points where $|\theta|^2$ and
$|\phi_1|$ vanish respectively, and we denote by $Z_{|\theta|^2}^c$ and $Z_{|\phi_1|}^c$
their corresponding complements, we will have that $Z_{|\phi_1|}^c\subset
Z_{|\theta|^2}$, therefore if $[A_1,\phi_1]$ is not a reducible solution then
$Z_{|\phi_1|}^c$ is a non-empty open set. By a result of N. Aronszajn (see \cite{aronszajn})
we will have that $\theta=0$, since it vanishes in an open set.
\end{proof}


Since $\mathcal{C}(\mathfrak{c})$ is an affine space it is contractible. Also the space
of reducible configurations $\mathcal{A}(\mathfrak{c})\times\{0\}$ is contractible and
has infinite codimension in $\mathcal{C}(\mathfrak{c})$. Since
$\mathcal{C}^*(\mathfrak{c})$ is open in $\mathcal{C}(\mathfrak{c})$ and it is the
complement of $\mathcal{A}(\mathfrak{c})\times\{0\}$ then it is contractible.
$\mathcal{B}^*(\mathfrak{c})=\mathcal{C}^*(\mathfrak{c})/\mathcal{G}(\mathfrak{c})$ is
the classifying space of $\mathcal{G}(\mathfrak{c})=Map(M,S^1)$ since
$\mathcal{G}(\mathfrak{c})$ acts freely on $\mathcal{C}^*(\mathfrak{c})$.

Moreover,
\begin{equation*}
Map(M,S^1)\sim Map(M,S^1)_o\times\pi_0(Map(M,S^1)),
\end{equation*}
where $Map(M,S^1)_o$
denotes homotopically constant maps. $Map(M,S^1)_o$ can be identified with $S^1$,
therefore $Map(M,S^1)\sim S^1\times H^1(M;\Z)$, so the classifying space for $Map(M,S^1)$
is weakly homotopically equivalent to $\CP{\infty}\times\frac{H^1(M;\R)}{H^1(M'\Z)}$, and
\begin{equation}\label{cmoduli}
H^*(\mathcal{B}^*(\mathfrak{c});\Z)\cong\Z[U]\otimes\Omega^*H^1(M;\Z),
\end{equation}
where $U$ is a generator for $H^*(\CP{\infty};\Z)$.

\begin{definition}
The \mbox{Seiberg-Witten} invariant $SW(\mathfrak{c})$ for the $\spin^c$-structure
$\mathfrak{c}$ is defined as follows
\begin{equation*}
SW(\mathfrak{c})=\begin{cases} \langle
U^{d(\mathfrak{c})/2},\mathcal{M}(\mathfrak{c}\rangle|_{\mathcal{B}^*(\mathfrak{c})}
&\quad\text{if }d(\mathfrak{c})\text{ is even}\\ 0 &\quad\text{otherwise}
\end{cases}
\end{equation*}
\end{definition}
It is easy to see that this invariant is a cobordism invariant of the moduli space
$\mathcal{M}(\mathfrak{c})$, therefore it does not depend on the metric we used to define
the Dirac operator, it does define an invariant of the smooth manifold $M$.

From this definition it is easy to see that we are loosing information about the moduli
space. For example if the moduli space is odd dimensional this invariant is zero, even
though the moduli itself may not represent a trivial bordism class in
$\mathcal{B}^*(\mathfrak{c})$.

\begin{definition}\label{monopole}
Let $(M,\mathfrak{c})$ be a smooth compact oriented \mfld{4} with a $\spin^c$-structure
$\mathfrak{c}$. We will say that $\mathfrak{c}$ is a \emph{$B$-class} if for some (then
for any) Riemannian metric $g$ on $M$, the moduli space $\mathcal{M}_g(\mathfrak{c})$ of
irreducible solutions to the SW-equations is a smooth manifold of dimension
$d(\mathfrak{c})\ge 0$ that represents a non-trivial bordism class in
$\mathcal{B}^*(\mathfrak{c})$, {\it i.e.} there exists $\eta\in
H^*(\mathcal{B}^*(\mathfrak{c});\Z)$ of degree $d(\mathfrak{c})$ such that
\begin{equation*}
\langle\eta,\mathcal{M}(\mathfrak{c})\rangle|_{\mathcal{B}^*(\mathfrak{c})}\neq 0.
\end{equation*}
\end{definition}

\section{SW-Equations and Conformal Structures}

It is easy to see that conformal changes on the metric can be lifted to a fixed
$\spin^c$-structure, and one can study the associated change in the Dirac operator. A
basic important fact is that \emph{the Dirac operator remains essentially invariant under
all conformal changes of the metric}.

We now make this statement precise. Let $(M,\mathfrak{c})$ be a
fixed smooth compact oriented \mfld{n} with a fixed
$\spin^c$-structure $\mathfrak{c}$ and a fixed Hermitian structure
$h$ on the determinant line bundle $L_\mathfrak{c}$. Fix a
Riemannian metric $g$ on $M$ and consider the conformally related
metric $g_f=e^{2f}g$, where $f$ is a smooth function on $M$.  To
each $g$-orthonormal tangent frame $\{e_i\}_{i=1\ldots n}$ we can
associate the $g_f$-orthonormal frame $\{e'_i\}_{i=1\ldots n}$,
where $e'_i=\psi_f(e_i)=e^{-f}e_i$ for each $i$. This map induces
a bundle isometry between the bundles $S(\mathfrak{c})$ and
$S'(\mathfrak{c})$. Let $\Psi_f=e^{-\frac{n-1}{2}f}\psi_f$. The
resulting map is a bundle isomorphism which is conformal on each
fiber. The proof of the following proposition is similar to the
one found in \cite{law&mich}, pages $132-134$, since we are not
changing the $U(1)$-connection on $L_\mathfrak{c}$.

\begin{proposition}\label{cDirac}
Let $D_A$ and $D'_A$ be the Dirac operators (induced by the
$U(1)$-connection $A$) defined over the conformally related
Riemannian manifolds \linebreak $(M,g)$ and $(M,g_f)$
respectively. Then
\begin{equation*}
\Psi_f\circ D_A=D'_A\circ\Psi_f
\end{equation*}
\end{proposition}

\begin{corollary}
There is bijection between $\ker D_A$ and $\ker D'_A$.
\end{corollary}

Let $(M,\mathfrak{c})$ be a fixed smooth compact oriented \mfld{4}
with a fixed $\spin^c$-structure $\mathfrak{c}$. We want to relate
the moduli spaces $\mathcal{M}(\mathfrak{c})$ and
$\mathcal{M}'(\mathfrak{c})$ for two Riemannian metrics $g$ and
$g_f$ (respectively) in the same conformal class. It is well known
(see \cite{morgan:sw}) that both moduli spaces represent the same
bordism class (in $\mathcal{B}^*(\mathfrak{c})$), but when one of
the metrics is K\"ahler, both moduli spaces are diffeomorphic (see
proposition \ref{difmoduli}) .

\begin{proposition}\label{1to1}
Let $(M,\mathfrak{c})$ be a fixed smooth compact oriented \mfld{4} with a fixed
$\spin^c$-structure $\mathfrak{c}$. Let $g$ be a fixed Riemannian metric on $M$ and
consider the conformal metric $g_f=e^{2f}g$. Solutions to the \mbox{Seiberg-Witten}
equation for the metric $g_f$ are in one-to-one correspondence with solutions of the
following pair of equations:
\begin{equation}\label{confsw}
\begin{split}
D_A\phi&=0\\ F_A^+&=e^{-f}q(\phi).
\end{split}\tag{$SW_f$}
\end{equation}
The one-to-one correspondence is given by the map $(A,\phi)\mapsto (A,\Psi_f\phi)$.
\end{proposition}
\begin{proof}
This is a consequence of Proposition \ref{cDirac}, the expression for $q$ (see Definition
\ref{defsw}) and that $\star'|_{\bigwedge^2}=\star|_{\bigwedge^2}$, where $\star$ and
$\star'$ are the Hodge operators of $g$ and $g_f$, respectively.
\end{proof}

\begin{proposition}\label{difmoduli}
Let $(M,g)$ be a K\"ahler surface with K\"ahler metric $g$. Then
for any smooth function $f:M\to\R$
\begin{itemize}
\item If the degree of $K_M$ is negative the only solutions to (\ref{confsw}) are reducible, \textit{i.e. } $\mathcal{M}_{e^{2f}g}(\mathfrak{c})=\emptyset$.
\item Let $\mathfrak{c}$ be the $\spin^c$-structure determined by the complex structure. If the degree of $K_M$ is positive then
$\#\mathcal{M}_{e^{2f}g}(\mathfrak{c})=1$.
\end{itemize}
\end{proposition}
\begin{proof}
The proof of this proposition can be carry out following the steps
in the proof of \textbf{Proposition 7.3.1} in \cite{morgan:sw} pg.
119, replacing the expression for $q$ with $e^{-f}q$.
\end{proof}

\begin{remark}
Note that $\#\mathcal{M}_{e^{2f}g}(\mathfrak{c})=1$ is stronger
than $SW_{e^{2f}g}(\mathfrak{c})=1$, which we already knew (see
\cite{morgan:sw}).
\end{remark}

\section{SW-Moduli Space of a Manifold with a
Cylindrical End}

The last result shows that if $(M,g)$ is a K\"ahler surface with
$\deg(K_M)>0$ the \mbox{Seiberg-Witten} moduli space for any
metric $g_f=e^{2f}g$ in the same conformal class of $g$ consists
of a single point. In this Section we extend this result to a
manifold with finitely many cylindrical ends.
\begin{definition}
We will say that $(M_\infty,g_\infty)$ is a \emph{manifold with a cylindrical end modeled
on $\R^+\times S^3$}, if $M_\infty$ is diffeomorphic to $M-\{p\}$ where $M$ is a closed
manifold, and $F:U_p-\{p\}\to\R^+\times S^3$ where $F(x)=(\log(|x|^{-1}),x/|x|)$ is a
diffeomorphism such that $(g_\infty)|_{U_p-\{p\}}$ is the $F$-pull-back of the standard
product metric $dt^2+g_{S^3}$ on $\R^+\times S^3$ and $U_p$ is a neighborhood of $p$.
\end{definition}

If $(M,g)$ is a Riemannian manifold such that $g$ is flat in a $\delta$-neighborhood of
$p$, where $\delta<\inj(M,g)$, there is a canonical way to produce a manifold with a
cylindrical end using the conformal class of $g$. Here $\inj(M,g)$ denotes the
injectivity radius of $(M,g)$. Choose a function $\lambda_l:(0,1]\to[1,\infty)$ which satisfies
\begin{equation}\label{funcl}
\lambda_l(r)=\begin{cases} 1 & \text{if } 0\le r\le e^{-l}\delta^3\\ \delta^2/r &\text{if }
e^{-l}\delta^2\le r\le\delta^2\\ 1 &\text{if } r\ge\delta.
\end{cases}
\end{equation}
Consider the sequence of functions $\{f_l\}$, where
$e^{f_l(x)}=\lambda_l(|x|)$ and the sequence of metrics
$g_l=e^{2f_l}g$. This sequence of metrics converges in the
$\mathcal{CO}$-topology on $M-\{p\}$ to a metric $g_\infty$. The
pair $(M-\{p\},g_\infty)$ is a manifold with a cylindrical end. We
will denote by $\Psi_l$ the associated conformal isomorphism
defined above proposition \ref{cDirac}.

The SW-equations make perfectly good sense on a manifold with a cylindrical end, but in
order to use the usual analytical tools, one has to extend the $\mathcal{C}^\infty$
objects to appropriate weighted Sobolev spaces (see \cite{mcowen}). From now on every
time we work on a manifold with finitely many cylindrical ends we extend the
configuration space $\mathcal{A}(\mathfrak{c})$ and the gauge group
$\mathcal{G}(\mathfrak{c})$ by requiring $A$ and $\phi$ to be in
$L^2_{2,\epsilon}(M_\infty,g_\infty)$ and $\sigma$ to be in
$L^2_{3,\epsilon}(M_\infty,g_\infty)$. The $L^p_{q,\epsilon}(M_\infty,g_\infty)$ norm is
defined as
\begin{equation*}
  \|h\|_{p,q,\epsilon}=\|e^{\tilde{\epsilon}t}h\|_{p,q},
\end{equation*}
where $\tilde{\epsilon}$ is a smooth non-decreasing function with bounded derivatives,
$\tilde{\epsilon}:M\to[0,\epsilon]$, such that $\tilde{\epsilon}(x)\equiv 0$ for $x\notin
B_\delta(p)$ and $\tilde{\epsilon}(x)\equiv \epsilon>0$ for $x\in B_{\delta^2}(p)$.

Here we choose the weight $\epsilon<1$ because we want to produce
solutions on the manifold with cylindrical end from solutions on
the manifold $(M,g)$ via the conformal process ($g_l\to g_\infty$)
using proposition \ref{1to1}.

\begin{proposition}
Let $(M,g)$ be any Riemannian \mfld{4}, where $g$ is flat in the
neighborhood of some point $p\in M$ . If $(A,\phi)$ is a solution
of (\ref{confsw}) on $(M,g)$ (where $f=f_\infty$) then
$(A,\Psi_\infty\phi)$ is a solution of the SW-equations on
$(M_\infty, g_\infty)$, such that $(A,\Psi_\infty\phi)\in
L^2_{1,\epsilon}(M_\infty,g_\infty)$.
\end{proposition}
\begin{proof}
The fact that $(A,\Psi_\infty\phi)$ satisfies the SW-equations
follows from proposition \ref{1to1}. We just need to show that
$(A,\Psi_\infty\phi)\in L^2_{1,\epsilon}(M_\infty,g_\infty)$. In
order to do this, we will use the metric $g$ as the background
metric.
\begin{align*}
  \|\Psi_\infty\phi\|^2_{2,1,\epsilon}&=\|e^{\tilde{\epsilon}t}\Psi_\infty\phi\|^2_{2,1}\\
  &\quad\int_{M-B_\delta(p)}(|\phi|^2+|\nabla\phi|^2)d\mu+\\
  &\quad\int_{\R^+\times S^3}(|e^{\tilde{\epsilon}t}\Psi_\infty\phi|^2_\infty+|e^{\tilde{\epsilon}t}\partial_t\nabla\Psi_\infty\phi|^2_\infty)dtd\mu_{S^3} \\
  &=\int_{M-B_\delta(p)}(|\phi|^2+|\nabla\phi|^2)d\mu+\\
  &\quad\int_{B_\delta(p)-\{p\}}(|r^{-\epsilon+3/2}\phi|^2+|r^{-\epsilon+1+3/2}\partial_r\nabla\phi|^2)\frac{1}{r}drd\mu_{S^3}\\
  &=\int_{M-B_\delta(p)}(|\phi|^2+|\nabla\phi|^2)d\mu+\\
  &\quad\int_{B_\delta(p)-\{p\}}r^{-2\epsilon-1}(|\phi|^2+|r\partial_r\nabla\phi|^2)r^3drd\mu_{S^3}\\
  &\leq C\|\phi\|^2_{2,1}.
\end{align*}
To prove that $A\in L^2_{1,\epsilon}(M_\infty,g_\infty)$ we need to recall that
\begin{equation*}
  \star_{g_f}|_{\wedge^p}=e^{(n-2p)f}\star_g|_{\wedge^p}
\end{equation*}
where $g_f=e^{2f}g$. The computation is very similar to the one above.
\end{proof}

Our next task is to show that there is no loss of generality in assuming that a K\"ahler metric $g$ is
flat in a neighborhood of some point.

\begin{proposition}\label{non-obstruction}
Let $(M^{2n},g)$ be a K\"ahler \mfld{2n} with K\"ahler metric $g$ and induced K\"ahler
form $\omega$ . There is no local obstruction to finding a K\"ahler metric on $M$, flat
in a neighborhood of a point (a finite collection of points) without changing the
K\"ahler class of $\omega$.
\end{proposition}
\begin{proof}
Let $p\in M$. The existence of such metric is equivalent to
finding a neighborhood $U$ of $p$, and a K\"ahler form $\omega'$
in the same K\"ahler class of $\omega$, such that
$\omega'|_U=\omega_0=\sum_{i=1}^ndz^i\wedge d\overline{z^i}$. It
is well known that there exist an $\epsilon$-neighborhood $U_p$ of
$p$ and a function $f:U_p\to\R$ such that
$\omega|_{U_p}=i\partial\overline{\partial}(z\overline{z}+f(z))>0$,
where $|f(z)|\sim o(|z|^4)$ and $|z|$ denotes the distance (using
the K\"ahler metric $g$) on $U_p$ to $p$. Let
$\mathcal{K}^\infty(f)$ be the space of smooth functions on $M$
that satisfy
\begin{equation*}
\mathcal{K}^\infty(f)=\{h_{s,t}\in\mathcal{C}^\infty(M)|\,h(z)=-f(z)\text{
if }|z|<s,\,h(z)=0\text{ if }t<|z|\}
\end{equation*}
where $0<s<t\leq\epsilon$, depend on $h$.
Observe that if $f$ is zero we do not have anything to prove,
otherwise $0\not\in\mathcal{K}^\infty(f)$, but
$0\in\mathcal{K}^{3+\alpha}(f)$, where $\mathcal{K}^{3+\alpha}(f)$
denotes the completion of $\mathcal{K}^\infty(f)$ in the
$\mathcal{C}^{3+\alpha}$ topology. To see this consider the
one-parameter family of functions $h_k(z)=-\rho(k|z|)f(z)$, where
$\rho$ 
 is a smooth bump function such that
\begin{equation*}
\rho(r)=\begin{cases} 1&\text{if } 0<r<1/2\\ 0&\text{if }1/2<r<1.
\end{cases}
\end{equation*}
All these functions are in $\mathcal{K}^\infty(f)$ and satisfy
\begin{align*}
|h_k(z)|&\sim o(|z|^4)\\ |\nabla h_k(z)|&\sim o(|z|^3)\\
|\nabla^2h_k(z)|&\sim o(|z|^2)\\ |\nabla^3h_k(z)|&\sim o(|z|)\\
|\nabla^4h_k(z)|&\sim o(1).
\end{align*}
It is not difficult to see that $h_k\to 0$ in the
$\mathcal{C}^{3+\alpha}$ topology. It is important to recall that
the set $\mathcal{P}(\omega)$ of smooth functions $h$ such that
$\omega_h=\omega+i\partial\overline\partial h>0$, is open in the
$\mathcal{C}^\infty$ topology. This two facts allow us to find
$h_{s,t}\in\mathcal{K}^\infty(f)\bigcap\mathcal{P}(\omega)$,
$\mathcal{C}^{3+\alpha}$ close to $0$, such that
\begin{align*}
\omega_{h_{s,t}}&=\omega+i\partial\overline\partial h_{s,t}>0\\
&=\omega_0+i\partial\overline\partial (f +h_{s,t}),
\end{align*}
therefore we have
\begin{equation*}
\omega_{h_{s,t}}|_{B_s(p)}=\omega_0,
\end{equation*}
where $B_s(p)=\{z\in U_p|\,|z|<s\}$.
\end{proof}

\begin{corollary}\label{excylend}
For any compact oriented K\"ahler surface $(M,g)$ with canonical line bundle $K_M$ of positive degree, where $g$ is flat
in a neighborhood of some point, the \emph{induced} manifold with a cylindrical end
$(M_\infty,g_\infty)$ admits solutions to the SW-equations.
\end{corollary}

In order to prove that the \mbox{Seiberg-Witten} moduli space of a
manifold with a cylindrical end consists of only one point if
$\deg(K_M)>0$, we will need the following technical result.

\begin{proposition}\label{cylend2compact}
Let $(M_\infty,g_\infty)$ be a \mfld{4} with a cylindrical end. If
\begin{equation*}(A_\infty,\phi_\infty)\in\mathcal{C}^\infty\bigcap
L^2_{k,\epsilon}(M_\infty,g_\infty)
\end{equation*}
is a solution of the SW-equations on the manifold with cylindrical
end $(M_\infty,g_\infty)$, then
$(A_\infty,\Psi_\infty^{-1}\phi_\infty)$ extends to a smooth
solution of (\ref{confsw}) on $(M,g)$, replacing the strictly
positive function $e^{-f}$ by the non-negative function
\begin{equation*}
\lambda_\infty(x)=\begin{cases} |x|/\delta^2&\text{if }|x|<\delta^2\\ 1&\text{if } |x|>\delta
\end{cases}
\end{equation*}
\end{proposition}
\begin{proof}
It is easy to see that $(A_\infty,\Psi_\infty^{-1}\phi_\infty)\in
L^2(M,g)$, as it is to see that
$(A_\infty,\Psi_\infty\phi_\infty)$ is a solution of
(\ref{confsw}) with function $\lambda_\infty$ replacing $e^{-f}$.
The first equation in (\ref{confsw}) tell us that
$\Psi_\infty^{-1}\phi_\infty$ is a holomorphic section on
$M-\{p\}$. Using Hartog's Theorem we can extend this to a
holomorphic section on $M$. All the analysis done in proving
proposition \ref{difmoduli} can be carry out if we replace the
strictly positive function $e^{-f}$ in (\ref{confsw}) by a
non-negative function $\lambda_\infty$ whose zero set has measure
zero.
\end{proof}

\begin{corollary}\label{cylendcor}
Let $(M,g)$ be a compact oriented K\"ahler surface with canonical line bundle $K_M$ of positive degree, where $g$ is flat
in a neighborhood of some point. Then there exists a solution
$(A_\infty,\phi_\infty)\in\mathcal{C}^\infty\bigcap L^2_{k,\epsilon}(M_\infty,g_\infty)$
of the SW-equations on $(M_\infty,g_\infty)$. This solution is unique up to gauge
equivalence.
\end{corollary}
\begin{proof}
Since all the analysis done in proving proposition \ref{difmoduli}
can be carry out if we replace the strictly positive function
$e^{-f}$ in (\ref{confsw}) by a non-negative function
$\lambda_\infty$ whose zero set has measure zero, existence is a
consequence of corollary \ref{excylend} and uniqueness is obtained
using proposition \ref{cylend2compact} and proposition
\ref{difmoduli}.
\end{proof}

\ifpdf
\section{Holonomy, Connected Sums and SW-Invariants}
\else
\section{Holonomy, Connected Sums with $S^1\times S^3$ and SW-Invariants}
\fi

Consider the diffeomorphism
\begin{equation*}
F:\R^4-\{0\}\to\R\times S^3,\quad F(x)=\left(\log|x|,\frac{x}{|x|}\right).
\end{equation*}
It is easy to see that the pull-back of the standard product metric $g$ on $\R\times S^3$
under this diffeomorphism is given by
\begin{equation*}
F^*g(\xi,\eta)=\frac{1}{|x|^2}\langle\xi,\eta\rangle
\end{equation*}
for $|x|\le 1$. Fix $\delta>0$ and choose a function $\lambda_l:(0,1]\to[1,\infty)$ as in (\ref{funcl}) and consider the metric
\begin{equation*}
g_l(\xi,\eta)=\lambda_l(|x|)^2\langle\xi,\eta\rangle.
\end{equation*}
Note that for $e^{-l}\delta^2\le|x|\le\delta^2$ this metric agrees with the above
pull-back metric $F^*g$.


It is convenient to think of the connected sum $M\#(S^1\times S^3)$ as follows. Let $M$
be a smooth compact oriented \mfld{4}. Fix two points $p_1,p_2\in M$, and choose a metric
$g$ on $M$ which is flat in a $\delta$-neighborhood of $p_i$. For 
every $l\in\N$ consider the $e^{-l-1}\delta^2$-neighborhood of $p_i$ (with respect to
$g$) $B_{p_i}(e^{-l-1}\delta^2)$, and denote by $M_l$ the open subset of $M$ given by the
complement of $\overline{B_{p_1}(e^{-l-1}\delta^2)}\cup
\overline{B_{p_2}(e^{-l-1}\delta^2)}$. If we denote by
$T_i=T_i(e^{-l}\delta^2,e^{-l-1}\delta^2)$ the annulus centered at $p_i$ with radii
$e^{-l-1}\delta^2$ and $e^{-l}\delta^2$, it is easy to see that there exist a
diffeomorphism (orientation reversing) that takes $T_1$ into $T_2$ and if we define
$g_l=\lambda_l^2g$, such diffeomorphism becomes a $g_l$-isometry. Since we have observed
that $T_1$ and $T_2$ are $g_l$-isometric we can identify $T_1$ with $T_2$ , and call them
$T_l$, to obtain a Riemannian manifold $(M\#_l(S^1\times S^3),g_l)$. 
 This manifold is simply the manifold $M$ with two cylindrical ends of
length $l$ obtained by conformally rescaling the metric $g$ and identifying the annuli.
It is easy to see that such manifold is diffeomorphic to the connected sum $M\#(S^1\times
S^3)$.


Even though the process above described can be realized on any
smooth \mfld{4} the following results are only valid when $M$ is a
K\"ahler surface, because to prove them, we (strongly) use that on
a given conformal class of metrics, the moduli spaces of solutions
of the SW-equations for any two representatives are diffeomorphic,
and this was proved for K\"ahler surfaces (see proposition
\ref{difmoduli}).

Our next task is to explain how a $\spin^c$-structure on $M$ transforms into a
$\spin^c$-structure on $M\#(S^1\times S^3)$ under the process above described. The
following Proposition will be very useful to explain it.

\begin{proposition}
There is a canonical projection map $\pi:M\#(S^1\times S^3)\to M$. It has the following properties:
\begin{enumerate}
\item\label{p1} The induced maps in cohomology
\begin{equation*}
\pi^*:H^i(M;\F)\to H^i(M\#(S^1\times S^3);\F)
\end{equation*}
are injective. Here $\F=\Z_2$ or $\Z$. Moreover for $i=0,2,4$, $\pi^*$ is an
isomorphism.
\item\label{p2} $\pi^*(w_2(M))=w_2(M\#(S^1\times S^3))$.
\end{enumerate}
\end{proposition}
We will denote the $\spin^c$-structure obtained in the above
proposition by $\mathfrak{c}_{0,1}$. It is not difficult to show
that the formal dimension of the moduli space associated to
$\mathfrak{c}_{0,1}$ is $d(\mathfrak{c}_{0,1})=d(\mathfrak{c})+1$.

To explain the increment in the dimension above we need to recall the concept of
\emph{holonomy}. Let $P_G\to M$ be a principal $G$-bundle over $M$, with a connection
$A$. Let $x\in M$ and denote by $C(x)$ the loop space at $x$. For each $\gamma\in C(x)$
the parallel displacement along $\gamma$ is an isomorphism of the fiber $\approx G$ onto
itself and we will denote it by $\hol_\gamma(A)$. The set of all such isomorphisms forms
a group, the \emph{holonomy group of $A$ with reference point $x$}.

Once and for all for each $l>0$ we will choose $p_l\in T_1$, $q_l\in T_2$ and a path
$\Gamma_l:I\to M$ from $p_l$ to $q_l$ such that after identifying $T_1$ with $T_2$ we
obtain and embedding $\gamma_l:S^1\to M\#_l(S^1\times S^3)$. It is not difficult to
observe that for all $l>0$ $[\gamma_l]\nsim 0\in \pi_1(M\#(S^1\times S^3))$, and in fact
$\gamma_l$ represents the $S^1$ factor of the connected sum.

If $A$ is a $U(1)$-connection on the determinant line bundle
$L_\mathfrak{c}$, we can trivialize $L_\mathfrak{c}$ along
$\Gamma_l$ so that the parallel transport along $\Gamma_l$ induces
the identity from the fiber at $p_l$ to the fiber at $q_l$. When
we identify $T_1$ with $T_2$ we still have the extra degree of
freedom of how to identify the fiber at $p_l$ with the fiber at
$q_l$, and this is measured by $\hol_\gamma(A)$, where $A$ is the
\emph{glued} connection. If we change of gauge, $\hol_\gamma(A)$
remains unchanged because the structure group $U(1)$ is Abelian.
In this section we will prove that when $M$ is a K\"ahler surface
then every solution to the \mbox{Seiberg-Witten} equations for a
$\spin^c$-structure $\mathfrak{c}$, induces an $S^1$ family of
solutions to the SW-equations for the $\spin^c$-structure
$\mathfrak{c}_{0,1}$ on $M\#(S^1\times S^3)$.

We can \emph{glue} a solution $(A_\infty,\phi_\infty)$ of the SW-equations on $(M_\infty,g_\infty)$ to produce a solution
$(A_l,\phi_l)$ of the following set of equations on $(M\#_l(S^1\times S^3),g_l)$
\begin{align*}
D_{A_l}\phi_l&=\mu(A_l,\phi_l)=\mu_l\\ F_{A_l}^+-q(\phi_l)&=\nu(A_l,\phi_l)=\nu_l,
\end{align*}
where $(\mu_l,\nu_l)\in\mathcal{S}(\mathfrak{c})\times\Omega^2_+(M\#_l(S^1\times
S^3);i\R)$. It is not difficult to see that
\begin{align*}
(\mu_l,\nu_l)\in L^2_1(M\#_l(S^1\times S^3),g_l)\\ \lim_{l\to\infty}\|(\mu_l,\nu_l)\|_{2,1}=0,
\end{align*}

\begin{definition}
We will denote by
$\mathcal{M}_\theta(\mathfrak{c}_{0,1})\subset\mathcal{M}(\mathfrak{c}_{0,1})$ the
solution subspace of the SW-equations satisfying the extra condition
\begin{equation*}
\hol_\gamma(A)=\theta,
\end{equation*}
and by $SW_\theta(\mathfrak{c}_{0,1})$ the cobordism invariant
associated to this moduli space (coun\-ting solutions with
appropriate sign). Note that the condition $\hol_\gamma(A)=\theta$
reduces the dimension of the moduli space by one.
\end{definition}

\begin{proposition}\label{propempty}
Let $(M\#_l(S^1\times S^3),g_l)$ be the connected sum of $M$ with $S^1\times S^3$ with a
neck of length $l$. For every $\theta\in S^1$ and for every $l\gg 0$, there exists some
generic perturbation $\eta_l\in\Omega^2_+(M\#(S^1\times S^3);i\R)$ with
$\supp\eta_l\subset T_l$ such that $SW_{\theta,l}^{-1}(0,\eta_l)\neq\emptyset$, where
$SW_{\theta,l}(A,\phi)=(D_A\phi,F_A^+-q(\phi))$ and $\hol_\gamma(A)=\theta$.
\end{proposition}
\begin{proof}
Observe that the condition of $\eta_l$ having $\supp\eta_l\subset T_l$ is not much of a
restriction at all, because the space of such $2$-forms is open and the set of generic
perturbations  is dense (see \cite{morgan:sw}).

Suppose otherwise, there exists some $\theta\in S^1$ such that for every $l\gg 0$ we have
$SW_{\theta,l}^{-1}(0,\eta_l)=\emptyset$. This would imply that
$SW_\infty^{-1}(0,0)=\emptyset$ since we have seen (see Corollary \ref{cylendcor}) that
$(M\#_l(S^1\times S^3),g_l)\to (M_\infty,g_\infty)$, but this is a contradiction because
we have proven (see Corollary \ref{cylendcor}), that $SW_\infty^{-1}(0,0)\neq\emptyset$.
\end{proof}

\begin{definition}
We will say that $(\tilde{A}_l,\tilde{\phi}_l)$ on $M_\infty$, $\mathcal{C}^0$-extends a
solution $(A_l,\phi_l)$ of $SW_{\theta,l}(A,\phi)=(0,\eta_l)$ on $M\#_l(S^1\times S^3)$
if
\begin{align*}
(\tilde{A}_l,\tilde{\phi}_l)|_{M_l}&\equiv (A_l,\phi_l) \text{
and}\\
(\tilde{A}_l(t,x),\tilde{\phi}_l(t,x))&=(A_l(x),e^{-2\epsilon
t}\phi_l(x))\\&\quad\quad \text{ for } (t,x)\in [l,\infty)\times
S^3.
\end{align*}
\end{definition}
\begin{remark}
Note that $(\tilde{A}_l,\tilde{\phi}_l)\in
L^2_{0,\epsilon}(M_\infty,g_\infty)$. From now on we will fix a
$U(1)$-connection $A$ on $L_\mathfrak{c}$.
\end{remark}


\begin{lemma}\label{lemexist}
If for every $l\gg 0$ there exist two different irreducible solutions $[A_l^1,\phi_l^1]$
and $[A_l^2,\phi_l^2]$ of
\begin{align*}
D_A\phi&=0\\ F_A^+-q(\phi)&=\eta_l
\end{align*}
on $M\#_l(S^1\times S^3)$ for some generic perturbations $\eta_l$, then
\begin{equation*}
(C_l,\psi_l)=(\|\tilde{\phi}_l^1-\tilde{\phi}_l^2\|_{2,0,\epsilon}(\tilde{A}_l^1-\tilde{A}_l^2),\frac{1}{\|\tilde{\phi}_l^1-\tilde{\phi}_l^2\|_{2,0,\epsilon}}(\tilde{\phi}_l^1-\tilde{\phi}_l^2))
\end{equation*}
satisfies
\begin{align*}
(C_l,\psi_l)\to (C,\psi)&\in L^2_{1,\epsilon}(M_\infty,g_\infty)\\ \|\psi\|_{2,0,\epsilon}&=1,
\end{align*}
where $(\tilde{A}_l^i,\tilde{\phi}_l^i)$ $\mathcal{C}^0$-extends $(A_l^i,\phi_l^i)$ to
$(M_\infty,g_\infty)$ for $i=1,2$, and $(A_l^i,\phi_l^i)$ are the unique representatives
obtained by the gauge fixing condition $\delta(A_l^i-A)=0$.
\end{lemma}

\begin{lemma}\label{lemuniq}
The same hypothesis as before. If $(\tilde{A}_l^i,\tilde{\phi}_l^i)\to
(A_\infty,\phi_\infty)$ in the $L^2_{1,\epsilon}(M_\infty,g_\infty)$ topology, then we have
\begin{equation*}
\lim_{l\to\infty}\|D(SW_{\theta,l})_{(\tilde{A}^1_l,\tilde{\phi}^1_l)}(C_l,\psi_l)\|_{2,0,\epsilon}=0
\end{equation*}
\end{lemma}
\begin{remark}
The proof of the previous two lemmas is not difficult but
technical (a straightforward computation) so we will omit the details.
\end{remark}

\begin{proposition}
For every $\theta\in S^1$, $SW_\theta(\mathfrak{c}_{0,1})=1$.
\end{proposition}
\begin{proof}
Assume that $SW_\theta(\mathfrak{c}_{0,1})\neq \pm 1$. Proposition
\ref{propempty} implies, for $l\gg 0$ there exist (at least) two
different irreducible solutions $(A_l^i,\phi_l^i)$, $i=1,2$ on
$(M\#_l(S^1\times S^3),g_l)$. By lemma \ref{lemexist} and lemma
\ref{lemuniq} we would have an element of $\ker DSW_\infty$ at
$(A_\infty,\phi_\infty)$ the unique solution on
$(M_\infty,g_\infty)$, obtained in corollary \ref{cylendcor}. But
this is a contradiction since $(A_\infty,\phi_\infty)$ is a smooth
point. The same kind of argument shows that
$SW_\theta(\mathfrak{c}_{0,1})=1$ since
$SW_\infty(\mathfrak{c})=1$.
\end{proof}

\ifpdf
\section{Cohomology}
\else
\section{Cohomology of $\mathcal{B}^*(\mathfrak{c})$}
\fi

In this section we will build cohomology classes for $\mathcal{B}^*(\mathfrak{c})$ in
order to detect $B$-classes (see Definition \ref{monopole}). To describe the cohomology
of $\mathcal{B}^*(\mathfrak{c})$ we have to introduce the concept of \emph{universal
family of $SW$-connections} associated to a $\spin^c$ structure $\mathfrak{c}$,
parameterized by $\mathcal{B}^*(\mathfrak{c})$. A $SW$-connection is simply a pair
$(A,\phi)$, where $A$ is a $U(1)$-connection on $L_\mathfrak{c}$ and $0\neq\phi\in
S^+(\mathfrak{c})$.

A cohomology class $\beta\in H^i(\mathcal{B}^*(\mathfrak{c});\Z)$ can be thought of as a
homomorphism $\beta:H_i(\mathcal{B}^*(\mathfrak{c});\Z)\to\Z$, and the elements of
$H_i(\mathcal{B}^*(\mathfrak{c});\Z)$ can be thought of as homotopic classes of maps
$f:T\to\mathcal{B}^*(\mathfrak{c})$, where $T$ is a compact space. The maps
$f:T\to\mathcal{B}^*(\mathfrak{c})$ are naturally interpreted in terms of families of
$SW$-connections.

\begin{definition}
\emph{A family of $SW$-connections in a bundle $L_\mathfrak{c}\to M$ pa\-ra\-me\-tri\-zed
by a space $T$} is a bundle $L\to T\times M$ with the property that each slice
$L_t=L|_{\{t\}\times M}$ is isomorphic to $L_\mathfrak{c}$, together with a
$SW$-connection $(A_\phi)_t=(A_t,\phi_t)$ in $L_t$, forming a family
$A_\phi=\{(A_\phi)_t\}$.
\end{definition}

Let $p_2:\mathcal{C}^*(\mathfrak{c})\times M\to M$ be the projection onto the second
factor and let $\mathcal{L}_\mathfrak{c}\to\mathcal{C}^*(\mathfrak{c})\times M$ be the
pull-back line bundle, $p_2^*L_\mathfrak{c}$. Then $\mathcal{L}_\mathfrak{c}$ carries a
tautological family of $SW$-connections $A_\phi$, in which the $SW$-connection on the
slice $\mathcal{L}_\mathfrak{c}|_{\{(A,\phi)\}}$ over $\{(A,\phi)\}\times M$ is
$(p_2^*(A),p_2^*(\phi))$. The group $\mathcal{G}(\mathfrak{c})$ acts freely on
$\mathcal{C}^*(\mathfrak{c})\times M$ as well as on
$\mathcal{L}_\mathfrak{c}=\mathcal{C}^*(\mathfrak{c})\times L_\mathfrak{c}$, and there is
therefore a quotient bundle
\begin{align*}
\mathbb{L}_\mathfrak{c}&\to\mathcal{B}^*(\mathfrak{c})\times M\\
\mathbb{L}_\mathfrak{c}&=\mathcal{L}_\mathfrak{c}/\mathcal{G}(\mathfrak{c}).
\end{align*}
The family of $SW$-connections $A_\phi$ is preserved by $\mathcal{G}(\mathfrak{c})$, so
$\mathbb{L}_\mathfrak{c}$ carries an inherited family of $SW$-connections
$\mathbb{A}_\phi$. This is the \emph{universal family of $SW$-connections in
$L_\mathfrak{c}\to M$ parameterized by $\mathcal{B}^*(\mathfrak{c})$}.

If a family of $SW$-connections is parameterized by a space $T$ and carried by a bundle
$L\to T\times M$, there is an associated map $f:T\to\mathcal{B}^*(\mathfrak{c})$ given by
\begin{equation*}
f(t)=[A_t,\phi_t].
\end{equation*}
Conversely, given $f:T\to\mathcal{B}^*(\mathfrak{c})$ there is a corresponding pull-back
family of connections carried by $(f\times I)^*\mathbb{L}_\mathfrak{c}$. These two
constructions are inverses of one another: if $f$ is determined by the above equation,
then for each $t$ there is a \emph{unique} isomorphism $\psi_t$ between the
$SW$-connections in $L_t$ and $(f\times I)^*(\mathbb{L}_\mathfrak{c})_t$, and as $t$
varies these fit together to form an isomorphism $\psi:L\to (f\times
I)^*\mathbb{L}_\mathfrak{c}$ between these two families. (The uniqueness of $\psi_t$
results from the fact that $\mathcal{G}(\mathfrak{c})$ acts freely on
$\mathcal{C}^*(\mathfrak{c})$). Thus:

\begin{lemma}\label{swfamily}
The maps $f:T\to\mathcal{B}^*(\mathfrak{c})$ are in one-to-one correspondence with
families of $SW$-connections on $M$ parameterized by $T$, and this correspondence is
obtained by pulling back from the universal family
$(\mathbb{L}_\mathfrak{c},\mathbb{A}_\phi)$.
\end{lemma}

\begin{remark}
Let $\{\gamma_i\}$ be fixed representatives for the generators of the free part of
$H_1(M;\Z)$. If $f_1,f_2:T\to \mathcal{B}^*(\mathfrak{c})$ are homotopic, the
corresponding bundles $L_1$ and $L_2$ are isomorphic, and the corresponding holonomy maps
$h_1:T\to (S^1)^{b_1}$ and $h_2:T\to (S^1)^{b_1}$ are homotopic, where the holonomy map
is defined as $h_i(t)=(\hol_{\gamma_1}(f_i(t)),\ldots,\hol_{\gamma_{b_1}}(f_i(t)))$.
\end{remark}

There is a general construction which produces cohomology classes in
$\mathcal{B}^*(\mathfrak{c})$, using the slant-product pairing
\begin{equation*}
/:H^{d-i}(\mathcal{B}^*(\mathfrak{c});\Z)\times H_i(M;\Z)\to
H^i(\mathcal{B}^*(\mathfrak{c});\Z).
\end{equation*}

We have built over $\mathcal{B}^*(\mathfrak{c})\times M$ a line bundle
$\mathbb{L}_\mathfrak{c}$, so we can define a map
\begin{equation*}
\mu:H_i(X;\Z)\to H^{2-i}(\mathcal{B}^*(\mathfrak{c});\Z)
\end{equation*}
by
\begin{equation*}
\mu(\alpha)=c_1(\mathbb{L}_\mathfrak{c})/\alpha.
\end{equation*}

If $T$ is any $(2-i)$-cycle in $\mathcal{B}^*(\mathfrak{c})$, the class $\mu(\alpha)$ can
be evaluated on $T$ using the formula
\begin{equation*}
\langle\mu(\alpha),T\rangle_{\mathcal{B}^*(\mathfrak{c})}=\langle
c_1(\mathbb{L}_\mathfrak{c}),T\times\alpha\rangle_{\mathcal{B}^*(\mathfrak{c})\times M},
\end{equation*}
which expresses the fact that the slant product is the adjoint of the cross-product
homomorphism. Next we will describe another way to build cohomology classes.

\begin{definition}
A closed curve $\gamma:S^1\to M$ induces a \emph{holonomy map}
\begin{equation*}
\hol_\gamma:\mathcal{B}^*(\mathfrak{c})\to S^1
\end{equation*}
defined as the holonomy of the $SW$-connections $A_\phi$ along $\gamma$. The pull-back of
the canonical class $d\theta$ of $S^1$ defines a cohomology class on
$H^1(\mathcal{B}^*(\mathfrak{c});\Z)$ which we will call the \emph{holonomy class along
$\gamma$}.
\end{definition}

\begin{proposition}\label{dhol}
The cohomology groups of $\mathcal{B}^*(\mathfrak{c})$ are generated by the image of the
map $\mu:H_i(X;\Z)\to H^{2-i}(\mathcal{B}^*(\mathfrak{c});\Z)$. Moreover, given
$\gamma\in H_1(M;\Z)$, $\mu(\gamma)$ is the holonomy class along $\gamma$,
$\hol_\gamma^*(d\theta)$.
\end{proposition}
\begin{proof}
First we will prove that if $\{\gamma_i\}$ are fixed representatives for the generators
for the free part of $H_1(M;\Z)$ then $\{\mu(\gamma_i)\}$ generates
$H^1(\mathcal{B}^*(\mathfrak{c});\Z)$. It is enough to prove that for every $i$ we can
find $\beta_i:S^1\to\mathcal{B}^*(\mathfrak{c})$ such that
$\langle\mu(\gamma_i),\beta_i\rangle|_{\mathcal{B}^*(\mathfrak{c})}=1$. Consider the line
bundle $\gamma_i^*L_\mathfrak{c}\to S^1$, and observe that there is no obstruction to
extend it to a line bundle $L\to S^1\times S^1$ such that $\deg L=\langle
c_1(L),S^1\times S^1\rangle=1$. Let $A_i$ be a $U(1)$-connection on $L$ and consider the
map
\begin{equation*}
\hol_{\,\bullet\times S^1}(A_i):S^1\to S^1.
\end{equation*}
It is not difficult to see that $\deg L=\deg(\hol_{\,\bullet\times
S^1}(A_i))$. After extending $A_i(t,\gamma_i)$ to a
$U(1)$-connection on $L_\mathfrak{c}\to M$ for each $t$, we obtain
(see remark below lemma \ref{swfamily}) our desired maps
$\beta_i:S^1\to\mathcal{B}^*(\mathfrak{c})$.

To prove the last statement we proceed as follows: let
$\alpha:S^1\to\mathcal{B}^*(\mathfrak{c})$,
\begin{align*}
\langle\mu(\gamma_i),\alpha\rangle_{\mathcal{B}^*(\mathfrak{c})}&=\langle
c_1(\mathbb{L}_\mathfrak{c}),\alpha\times\gamma_i\rangle_{\mathcal{B}^*(\mathfrak{c})\times
M}\\ &=\langle c_1((\alpha\times\gamma_i)^*(\mathbb{L}_\mathfrak{c})),S^1\times
S^1\rangle\\ &=\deg(\hol_{\,\bullet\times S^1}(A_i):S^1\to S^1)\\
&=\deg(\hol_{\gamma_i}\circ\alpha:S^1\to S^1)\\
&=\langle\deg_{\beta_i}^*(d\theta),\alpha\rangle_{\mathcal{B}^*(\mathfrak{c})}.
\end{align*}

Finally we have to show that if $x\in M$ then $\mu(x)$ generates the cohomology of the
$\CP{\infty}$ factor. Since $Map(M,S^1)_o$ acts freely on $\mathcal{C}^*(\mathfrak{c})$,
then it is easy to show that
$\mathbb{L}_\mathfrak{c}|_{\mathcal{B}^*(\mathfrak{c})}\approx\mathcal{C}^*(\mathfrak{c})/\mathcal{G}_0(\mathfrak{c})$,
where $\mathcal{G}_0(\mathfrak{c})$ is the kernel of the homomorphism
$\mathcal{G}(\mathfrak{c})\to S^1$ given by evaluating on the fiber over $x$.
\end{proof}

\section{Applications}

C. LeBrun \cite{lebrun} showed that under some mild conditions on $M$,
$M\#k\,\overline{\CP{2}}$ does not admit \mbox{Einstein} metrics. The precise statement
is the following:

\begin{lebrun} Let $M$ be a smooth compact oriented \mfld{4} with \mbox{$2e+3\sigma>0$}. Assume, moreover, that $M$ has a non-trivial Seiberg-Witten invariant. If \mbox{$k\ge\frac{25}{57}(2e+3\sigma)$} then $M\#k\,\overline{\CP{2}}$ does not admit an \mbox{Einstein} metric.
\end{lebrun}

\begin{remark}
The proof of this theorem only requires that $M$ has a
$\spin^c$-structure $\mathfrak{c}$ that is a $B$-class.
\end{remark}

\begin{theoremA}\label{gluing}
Let $(M,\mathfrak{c})$ be a smooth compact K\"ahler surface with a $\spin^c$-structure
$\mathfrak{c}$. There is a canonical $\spin^c$ structure in the connected sum manifold
$M\#(S^1\times S^3)$ which we will denote by $\mathfrak{c}_{0,1}$. Moreover
$d(\mathfrak{c}_{0,1})=d(\mathfrak{c})+1$. If $\mathfrak{c}$ is a non-trivial SW-class
for $M$ then $\mathfrak{c}_{0,1}$ is a $B$-class for the connected sum $M\#(S^1\times
S^3)$.
\end{theoremA}
\begin{proof}
$SW_\theta(\mathfrak{c}_{0,1})$ is a cobordism invariant for every
$\theta\in S^1$. Consider the smooth cobordism induced by the
family of metrics $g_l$ on $M\#(S^1\times S^3)$ as $l\to\infty$
and observe (corollary \ref{cylendcor}) that
$SW_\infty(\mathfrak{c})=1$. This shows that
\begin{equation*}
\langle\hol_\gamma^*(d\theta),\mathcal{M}(\mathfrak{c}_{0,1})\rangle|_{\mathcal{B}^*(\mathfrak{c}_{0,1})}=1,
\end{equation*}
where $\gamma$ is a representative for the $S^1$ factor of the
connected sum. This, the definition of a $B$-class and Proposition
\ref{dhol} complete the proof.
\end{proof}

\begin{corollary}
Let $(M,\mathfrak{c})$ be a smooth compact oriented K\"ahler surface with a
$\spin^c$-structure $\mathfrak{c}$. There is a canonical $\spin^c$ structure in the
connected sum $M\#2(S^1\times S^3)$ which we will denote by $\mathfrak{c}_{0,2}$.
Moreover $d(\mathfrak{c}_{0,2})=d(\mathfrak{c})+2$. If $\mathfrak{c}$ is a non-trivial
SW-class then $\mathfrak{c}_{0,2}$ is a $B$-class but has trivial \mbox{Seiberg-Witten}
invariant.
\end{corollary}
\begin{proof}
Theorem A shows that every time that we perform a connected sum with
$S^1\times S^3$ we \emph{add a cycle} to the moduli space, that lies entirely in the
$H^1(M\#(S^1\times S^3);\R)/H^1(M\#(S^1\times S^3);\Z)$ part of
$\mathcal{B}^*(\mathfrak{c}_{0,1})$.
\end{proof}

\begin{lemma}\label{lem}
Let $(M,\mathfrak{c})$  be a smooth compact oriented K\"ahler surface with a
$\spin^c$-structure $\mathfrak{c}$ and \mbox{$2e+3\sigma>0$}. Assume that $\mathfrak{c}$
is a non-trivial SW-class. Let $k,l$ be any two natural numbers. Then there is a
$B$-class $\mathfrak{c}_{k,l}$ on $M_{k,l}=M\#k\,\overline{\CP{2}}\#\,l(S^1\times S^3)$
such that
\begin{equation*}
(c_1^+(\mathfrak{c}_{k,l}))^2\ge(2e+3\sigma)(M).
\end{equation*}
\end{lemma}
\begin{proof}
First observe that $M_{k,l}=(M\#k\,\overline{\CP{2}})_{0,l}$. Since $M$ is a K\"ahler
surface, we know that $M\#k\,\overline{\CP{2}}$ is also a K\"ahler surface, and its
associated $\spin^c$ structure $\mathfrak{c}_{k,0}$ satisfies
$c_1(\mathfrak{c}_{k,0})=c_1(\mathfrak{c})+\sum_{j=1}^kE_j$, where $E_1,\ldots,E_k$ are
generators for the pull-backs to $M\#k\,\overline{\CP{2}}$ of the $k$ copies of
$H^2(\CP{2},\Z)$ so that
\begin{equation*}
c_1^+(\mathfrak{c})\cdot E_j\ge 0,\quad j=1,\ldots,k.
\end{equation*}
Let $c_1(\mathfrak{c}_{k,l})$ be the first Chern class of
$(\mathfrak{c}_{k,0})_{0,l}$ which is a $B$-class by theorem A, and notice that
$c_1(\mathfrak{c}_{k,l})=c_1(\mathfrak{c}_{k,0})$. One then has
\begin{align*}
(c_1^+(\mathfrak{c}_{k,l}))^2&=(c_1^+(\mathfrak{c}_{k,0}))^2\\
&=\left(c_1^+(\mathfrak{c})+\sum_{j=1}^kE_j^+\right)^2\\
&=(c_1^+(\mathfrak{c}))^2+2\sum_{j=1}^kc_1^+(\mathfrak{c}_{0,l})\cdot
E_j^++(\sum_{j=1}^kE_j^+)^2\\ &\ge (c_1^+(\mathfrak{c}))^2\\ &\ge (c_1(\mathfrak{c}))^2\\
&=(2e+3\sigma)(M).
\end{align*}
\end{proof}
LeBrun's theorem can be generalized in the following way:

\begin{theorem}\label{lebrungen}
Let $(M,\mathfrak{c})$ be a smooth compact oriented K\"ahler surface with a
$\spin^c$-structure $\mathfrak{c}$ and \mbox{$2e+3\sigma>0$}. Assume  that $\mathfrak{c}$
is a $B$-class. If $k+4l\ge\frac{25}{57}(2e+3\sigma)$ then
$M_{k,l}=M\#k\,\overline{\CP{2}}\#\,l(S^1\times S^3)$ does not admit an \mbox{Einstein}
metric.
\end{theorem}
\begin{proof}
The proof is the same as the one given by C. LeBrun \cite{lebrun}.
\end{proof}
There exists two well known topological obstructions to the
existence of Einstein metrics on a differentiable compact oriented
\mfld{4} $M$.

The first one is Thorpe's inequality (see \cite{besse2}), that comes from the
Gauss-Bonnet-Chern formula for the Euler characteristic $e(M)$ of
$M$ and from the Hirzebruch formula for the signature $\sigma(M)$
of $M$, which allow us to express these two
topological invariants in terms of the irreducible components of
the curvature under the action of $SO(4)$. It can be stated in the
following way
\begin{hitchin-thorpe}
Let $M$ be a compact oriented manifold of dimension $4$. If
$e(M)<\frac{3}{2}|\sigma(M)|$ then $M$ does not admit any Einstein
metric. Moreover, if $e(M)=\frac{3}{2}|\sigma(M)$ then $M$ admits
no Einstein metric unless it is either flat, or a $K3$ surface, or
an Enriques surface, or the quotient of an Enriques surface by a
free antiholomorphic involution.
\end{hitchin-thorpe}
This theorem implies a previous result of M. Berger who proved
that there exists no compact Einstein \mfld{4} with a negative
Euler characteristic. On the other hand, combining the
Gauss-Bonnet-Chern formula for the Euler characteristic with
Gromov's estimation of simplicial volume $\|M\|$ of a Riemannian
manifold $M$ (see \cite{gromov}), M. Gromov obtained the following
obstruction

\begin{gromov}
Let $M$ be a compact manifold of dimension $4$. If
$e(M)<\frac{1}{2592\pi^2}\|M\|$ then $M$ does not admit any
Einstein metric.
\end{gromov}


A. Sambusetti (see \cite{sambusetti}) found a topological obstruction to the existence of
Einstein metrics on compact $4$-manifolds which admit a non-zero degree map onto some
compact real or complex hyperbolic \mfld{4}. As a consequence, by connected sums, he
produces infinitely many non-homeomorphic $4$-manifolds which admit no Einstein metrics. This fact is not a consequence
of Hitchin-Thorpe's or Gromov's obstruction theorems. A. Sambusetti also
proves that any Euler characteristic and signature can be simultaneously realized by
these non-homeomorphic manifolds admitting no Einstein metrics.

\begin{definition}
We say that a pair $(m,n)\in\Z^2$ is \emph{admissible} if there
exists a smooth compact oriented \mfld{4} with Euler
characteristic $m$ and signature $n$. In fact a necessary and
sufficient condition for $(m,n)\in\Z^2$ to be an admissible pair
is that $m\equiv n\mod 2$.
\end{definition}
\noindent To prove our last result we need the following theorem
by Z. Chen (see \cite{chen}).

\begin{chen}
Let $x$, $y$ be integers satisfying
\begin{align*}
\frac{352}{89}x+140.2x^{2/3}<y&<\frac{18644}{2129}x-365.7x^{2/3},\\
x&>C,
\end{align*}
where $C$ is a large constant. There exists a simply connected
minimal surface $M$ of general type with $c_1^2(M)=y$,
$\chi(M)=x$. Furthermore, $M$ can be represented by a surface
admitting a hyperelliptic fibration.
\end{chen}

\begin{remark}
Recall that $\chi(M)$ denotes the Euler-Poincar\'e characteristic of the invertible sheaf $\mathcal{O}_M$. Using Noether's formula we have that
\begin{align*}
\chi(M)&=\frac{c_1^2(M)+e(M)}{12}\\
&=\frac{e(M)+\sigma(M)}{4}.
\end{align*}
If $M$ is not a complex surface $e(M)+\sigma(M)$ is not
necessarily a multiple of $4$ but it is always an even number.
\end{remark}

\begin{theoremB}
For each admissible pair $(m,n)$ there exist an infinite number of non-homeomorphic compact oriented $4$-manifolds which have Euler characteristic $m$ and signature $n$, with free fundamental group and which do not admit Einstein metric.
\end{theoremB}
\begin{proof}
Let $(m_0,n_0)$ be an admissible pair and consider the pair of
integers $(x'_0,y_0)=\left(\frac{m_0+n_0}{2},2m_0+3n_0\right)$. It
is always possible to find (infinitely many) positive integers $k$ and $l$ such that
\begin{align*}
  (x,y)=\left(\frac{x'_0+l}{2},y_0+k\right)&\in\mathcal{Z}\\
  4l+\frac{32}{57}k&\geq\frac{25}{57}y_0
\end{align*}
where $\mathcal{Z}$ denotes the set of $(x,y)\in\Z^2$ that satisfy
the conditions of Chen's theorem. The reason for this last
statement is that the region $\mathcal{Z}_\R$ determine by
$(x,y)\in\R^2$ such that
\begin{align*}
\frac{352}{89}x+140.2x^{2/3}<y&<\frac{18644}{2129}x-365.7x^{2/3},\\
x&>C,
\end{align*}
is open, connected and not bounded, where 
$C$ is the same constant as in Chen's theorem.

If we denote by $M$ the simply connected K\"ahler surface with
$c_1^2=y$ and $\chi=x$, then
$M_{k,l}=M\#k\,\overline{\CP{2}}\#l(S^1\times S^3)$, is a manifold
that realizes the pair $(m_0,n_0)$ and does not admit any Einstein
metric. This last statement is a consequence of theorem
\ref{lebrungen}.
\end{proof}


\end{document}